\documentclass[12pt]{article}
\usepackage{amssymb,amsmath}
\usepackage{hyperref}
\usepackage{graphicx}
\usepackage{color}
\usepackage{caption}

\begin{document}

\title{\LARGE\bf A set of the Vi\`ete-like recurrence relations for the unity constant}

\author{
\normalsize\bf S. M. Abrarov\footnote{\scriptsize{Dept. Earth and Space Science and Engineering, York University, Toronto, Canada, M3J 1P3.}}\, and B. M. Quine$^{*}$\footnote{\scriptsize{Dept. Physics and Astronomy, York University, Toronto, Canada, M3J 1P3.}}}

\date{February 3, 2017}
\maketitle

\begin{abstract}
Using a simple Vi\`ete-like formula for $\pi$ based on the nested radicals $a_k = \sqrt{2 + a_{k-1}}$ and $a_1 = \sqrt{2}$, we derive a set of the recurrence relations for the constant $1$. Computational test shows that application of this set of the Vi\`ete-like recurrence relations results in a rapid convergence to unity.
\vspace{0.25cm}
\end{abstract}

\section{Description and implementation}

\subsection{Derivation}

\vspace{0.25cm}
Several centuries ago the French mathematician Fran\c{c}ois Vi\`ete derived a remarkable formula for pi
\begin{equation}\label{eq_1}
\frac{2}{\pi }=\frac{\sqrt{2}}{2}\frac{\sqrt{2+\sqrt{2}}}{2}\frac{\sqrt{2+\sqrt{2+\sqrt{2}}}}{2}\cdots.
\end{equation}
Nowadays this well-known equation is commonly regarded as the Vi\`ete's formula for pi \cite{Herschfeld1935, Gearhart1990, Levin2005, Kreminski2008}. The uniqueness of this formula is due to nested radicals consisting of square roots of twos only. Defining these nested radicals as
$$
{{a}_{1}}=\sqrt{2},
$$	
$$
{{a}_{2}}=\sqrt{2+\sqrt{2}},
$$	
$$
{{a}_{3}}=\sqrt{2+\sqrt{2+\sqrt{2}}}
$$
$$
\vdots  \\ 
$$
\[
{{a}_{k}}=\underbrace{\sqrt{2+\sqrt{2+\sqrt{2+\cdots +\sqrt{2}}}}}_{k\,\,\text{square}\,\,\text{roots}}
\]
the Vi\`ete's formula \eqref{eq_1} for pi can be rewritten in a compact form as follows
$$
\frac{2}{\pi }=\underset{k\to \infty }{\mathop{\lim }}\,\prod\limits_{k=1}^{K}{\frac{{{a}_{k}}}{2}}.
$$

There is a simple Vi\`ete-like formula for pi that can be represented in form \cite{Abrarov2016}
\begin{equation}\label{eq_2}
\frac{\pi }{{{2}^{k+1}}}=\arctan \left( \frac{\sqrt{2-{{a}_{k-1}}}}{{{a}_{k}}} \right),	\qquad\qquad	k\ge 2,
\end{equation}
From this formula it follows that
\footnotesize
\begin{equation}\label{eq_3}
\begin{aligned}
\frac{\pi }{{{2}^{3}}}+\frac{\pi }{{{2}^{4}}}+\frac{\pi }{{{2}^{5}}}\cdots &= \\
&\hspace{-0.85cm}\arctan \left( \frac{\sqrt{2-{{a}_{1}}}}{{{a}_{2}}} \right)+\arctan \left( \frac{\sqrt{2-{{a}_{2}}}}{{{a}_{3}}} \right)+\arctan \left( \frac{\sqrt{2-{{a}_{3}}}}{{{a}_{4}}} \right)+\,\,\cdots
\end{aligned}
\end{equation}
\normalsize
and because of the decreasing geometric series
$$
\frac{1}{{{2}^{3}}}+\frac{1}{{{2}^{4}}}+\frac{1}{{{2}^{5}}}\cdots =\frac{1}{4}
$$
the equation \eqref{eq_3} can be expressed in a more simplified form
\begin{equation}\label{eq_4}
\frac{\pi }{4}=\underset{K\to \infty }{\mathop{\lim }}\,\sum\limits_{k=1}^{K}{\arctan \left( \frac{\sqrt{2-{{a}_{k}}}}{{{a}_{k+1}}} \right)}.
\end{equation}

It is more convenient for our purpose to represent the equation \eqref{eq_4} as
\footnotesize
\[
\begin{aligned}
\frac{\pi }{4}&=\arctan \left( \frac{\sqrt{2-\sqrt{2}}}{\sqrt{2+\sqrt{2}}} \right)+\arctan \left( \frac{\sqrt{2-\sqrt{2+\sqrt{2}}}}{\sqrt{2+\sqrt{2+\sqrt{2}}}} \right)\\
&+\arctan \left( \frac{\sqrt{2-\sqrt{2+\sqrt{2+\sqrt{2}}}}}{\sqrt{2+\sqrt{2+\sqrt{2+\sqrt{2}}}}} \right)+\,\,\cdots
\end{aligned}
\]
\normalsize
or
\[
\begin{aligned}
\frac{\pi }{4} &= \arctan \left( {{b}_{1}} \right)+\arctan \left( {{b}_{2}} \right)+\arctan \left( {{b}_{3}} \right)\,\,\cdots  \\ 
& =\underset{K\to \infty }{\mathop{\lim }}\,\sum\limits_{k=1}^{K}{\arctan \left( {{b}_{k}} \right),}  
\end{aligned}
\]
where the arguments of the arctangent functions can be found by using the recurrence relations
$$
{{b}_{k}}=\frac{\sqrt{2-{{a}_{k}}}}{{{a}_{k+1}}}
$$
and
$$
{{a}_{k}}=\sqrt{2+{{a}_{k-1}}},	\quad {{a}_{1}}=\sqrt{2}.
$$

Since
$$
\arctan \left( 1 \right)=\frac{\pi }{4}
$$
we can also write
\begin{equation}\label{eq_5}
\arctan \left( 1 \right)=\underset{K\to \infty }{\mathop{\lim }}\,\sum\limits_{k=1}^{K}{\arctan \left( {{b}_{k}} \right)}.
\end{equation}

The right side of the equation \eqref{eq_5} consists of the infinite summation terms of the arctangent functions. We may attempt to exclude the infinite sum using the identity
\begin{equation}\label{eq_6}
\arctan \left( x \right)+\arctan \left( y \right)=\arctan \left( \frac{x+y}{1-xy} \right)
\end{equation}
repeatedly. Specifically, we employ the following recurrence relations that just reflects the successive application of the identity \eqref{eq_6} above
$$
{{c}_{k}}=\frac{{{c}_{k-1}}+{{b}_{k}}}{1-{{c}_{k-1}}{{b}_{k}}}, \qquad {{c}_{1}}={{b}_{1}}.
$$
This enables us to rewrite the equation \eqref{eq_5} as
\begin{equation}\label{eq_7}
\arctan \left( 1 \right)=\arctan \left( {{c}_{k}} \right)+\underset{L\to \infty }{\mathop{\lim }}\,\sum\limits_{\ell =k+1}^{L}{\arctan \left( {{b}_{\ell }} \right)}.
\end{equation}

According to the Maclaurin expansion series
$$
\arctan \left( {{b}_{\ell }} \right)={{b}_{\ell }}-\frac{b_{\ell }^{3}}{3}+\frac{b_{\ell }^{5}}{5}-\frac{b_{\ell }^{7}}{7}+\cdots ={{b}_{\ell }}+O\left( b_{\ell }^{3} \right).
$$
Since at $\ell \to \infty $ the variable ${{b}_{\ell }}\to 0$ and, therefore, due to negligible $ O\left( b_{\ell }^{3} \right)$ we can simply replace it by $\arctan \left( {{b}_{\ell }} \right)$ and then use the equation \eqref{eq_2} in order to find a ratio of the limit
\begin{equation}\label{eq_8}
\underset{\ell \to \infty }{\mathop{\lim }}\,\frac{{{b}_{\ell +1}}}{{{b}_{\ell }}}=\underset{\ell \to \infty }{\mathop{\lim }}\,\frac{\arctan \left( {{b}_{\ell +1}} \right)}{\arctan \left( {{b}_{\ell }} \right)}=\underset{\ell \to \infty }{\mathop{\lim }}\,\frac{\pi /{{2}^{\ell +2}}}{\pi /{{2}^{\ell +1 }}}=\frac{1}{2}.
\end{equation}

Consider the following infinite sequence
\begin{equation}\label{eq_9}
\left\{ {{b}_{1}},{{b}_{2}},{{b}_{3}},\ldots ,{{b}_{\ell }},\ldots  \right\}.
\end{equation}
According to the limit \eqref{eq_8} the ratio ${{b}_{\ell +1}}/{{b}_{\ell }}$ tends to $1/2$ with increasing index $\ell $. Consequently, it is not difficult to see now that
$$
\frac{{{b}_{2}}}{{{b}_{1}}}<\frac{{{b}_{3}}}{{{b}_{2}}}<\frac{{{b}_{4}}}{{{b}_{3}}}<\cdots < \frac{{{b}_{\ell +1}}}{{{b}_{\ell }}} < \cdots <\frac{1}{2}.
$$
In fact, the tendency of the ratio ${{b}_{\ell +1}}/{{b}_{\ell }}$ towards $1/2$ with increasing index $\ell $ is very fast. In particular, when the index $\ell $ is large enough, say at $\ell >10$, the sequence \eqref{eq_9} behaves almost like a decreasing geometric progression where a common ratio is $1/2$.

Since the index $k$ in the equation \eqref{eq_7} can be taken arbitrarily large, we can rewrite it in form
\begin{equation}\label{eq_10}
\arctan \left( 1 \right)=\underset{k\to \infty }{\mathop{\lim }}\,\left[ \arctan \left( {{c}_{k}} \right)+\underset{L\to \infty }{\mathop{\lim }}\,\sum\limits_{\ell =k+1}^{L}{{{b}_{\ell }}} \right].
\end{equation}

Taking into account that the ratio ${{b}_{\ell +1}}/{{b}_{\ell }}$ tends to but never exceeds $1/2$, we can conclude that the damping rate in the sequence \eqref{eq_9} is faster than that of in a decreasing geometric progression 
$$
\left\{ {{b}_{1}},\frac{{{b}_{1}}}{2},\frac{{{b}_{1}}}{{{2}^{2}}},\frac{{{b}_{1}}}{{{2}^{3}}},\cdots \frac{{{b}_{1}}}{{{2}^{\ell }}}\cdots  \right\}
$$
with fixed common ratio $1/2$. This signifies that
$$
\sum\limits_{\ell =k+1}^{L}{{{b}_{\ell }}}<\sum\limits_{\ell =k+1}^{L}{\frac{{{b}_{1}}}{{{2}^{\ell -1}}}}, \qquad\qquad L > k > 0,
$$
and since the limit of the decreasing geometric series
$$
\underset{L\to \infty }{\mathop{\lim }}\,\sum\limits_{\ell =k+1}^{L}{\frac{{{b}_{1}}}{{{2}^{\ell -1}}}}\to 0, \qquad	 k\to \infty,
$$
we prove that
$$
\underset{L\to \infty }{\mathop{\lim }}\,\sum\limits_{\ell =k+1}^{L}{{{b}_{\ell }}}\to 0,	 \qquad k\to \infty.
$$
As a consequence, the equation \eqref{eq_10} can be further simplified as
\[
\arctan \left( 1 \right)=\underset{k\to \infty }{\mathop{\lim }}\,\arctan \left( {{c}_{k}} \right)\Leftrightarrow 1=\underset{k\to \infty }{\mathop{\lim }}\,{{c}_{k}}.
\]
Thus, we can infer that the constant $1$ can be approached successively by increment of the index $k$ in a set of the Vi\`ete-like recurrence relations
\begin{equation}\label{eq_11}
\left\{ \begin{aligned}
  & {{a}_{1}}=\sqrt{2}, \\ 
 & {{a}_{k}}=\sqrt{2+{{a}_{k-1}}}, \\ 
 & {{b}_{k}}=\frac{\sqrt{2-{{a}_{k}}}}{{{a}_{k+1}}}, \\ 
 & {{c}_{1}}={{b}_{1}}, \\ 
 & {{c}_{k}}=\frac{{{c}_{k-1}}+{{b}_{k}}}{1-{{c}_{k-1}}{{b}_{k}}}, \\ 
\end{aligned} \right.
\end{equation}
such that ${{c}_{k\to \infty }}\to 1.$

\subsection{Computation}

Consider the first three elements from the sequence \eqref{eq_9}
$$
{{b}_{1}}=\frac{\sqrt{2-{{a}_{1}}}}{{{a}_{2}}}=\frac{\sqrt{2-\sqrt{2}}}{\sqrt{2+\sqrt{2}}},
$$		
$$
{{b}_{2}}=\frac{\sqrt{2-{{a}_{2}}}}{{{a}_{3}}}=\frac{\sqrt{2-\sqrt{2+\sqrt{2}}}}{\sqrt{2+\sqrt{2+\sqrt{2}}}}
$$
and	
$$
{{b}_{3}}=\frac{\sqrt{2-{{a}_{3}}}}{{{a}_{4}}}=\frac{\sqrt{2-\sqrt{2+\sqrt{2+\sqrt{2}}}}}{\sqrt{2+\sqrt{2+\sqrt{2+\sqrt{2}}}}}.
$$
Consequently, the corresponding first three values of the variable ${{c}_{k}}$ are
$$
{{c}_{1}}={{b}_{1}}=\frac{\sqrt{2-\sqrt{2}}}{\sqrt{2+\sqrt{2}}}=\text{0}\text{.41421356237309504880}\ldots,
$$		
$$
{{c}_{2}}=\frac{{{c}_{1}}+{{b}_{2}}}{1-{{c}_{1}}{{b}_{2}}}=\frac{\frac{\sqrt{2-\sqrt{2}}}{\sqrt{2+\sqrt{2}}}+\frac{\sqrt{2-\sqrt{2+\sqrt{2}}}}{\sqrt{2+\sqrt{2+\sqrt{2}}}}}{1-\frac{\sqrt{2-\sqrt{2}}}{\sqrt{2+\sqrt{2}}}\frac{\sqrt{2-\sqrt{2+\sqrt{2}}}}{\sqrt{2+\sqrt{2+\sqrt{2}}}}}=\text{0}\text{.66817863791929891999}\ldots
$$
and
\[
\begin{aligned}
{{c}_{3}} &= \frac{{{c}_{2}}+{{b}_{3}}}{1-{{c}_{2}}{{b}_{3}}}=\frac{\frac{\frac{\sqrt{2-\sqrt{2}}}{\sqrt{2+\sqrt{2}}}+\frac{\sqrt{2-\sqrt{2+\sqrt{2}}}}{\sqrt{2+\sqrt{2+\sqrt{2}}}}}{1-\frac{\sqrt{2-\sqrt{2}}}{\sqrt{2+\sqrt{2}}}\frac{\sqrt{2-\sqrt{2+\sqrt{2}}}}{\sqrt{2+\sqrt{2+\sqrt{2}}}}}+\frac{\sqrt{2-\sqrt{2+\sqrt{2+\sqrt{2}}}}}{\sqrt{2+\sqrt{2+\sqrt{2+\sqrt{2}}}}}}{1-\frac{\frac{\sqrt{2-\sqrt{2}}}{\sqrt{2+\sqrt{2}}}+\frac{\sqrt{2-\sqrt{2+\sqrt{2}}}}{\sqrt{2+\sqrt{2+\sqrt{2}}}}}{1-\frac{\sqrt{2-\sqrt{2}}}{\sqrt{2+\sqrt{2}}}\frac{\sqrt{2-\sqrt{2+\sqrt{2}}}}{\sqrt{2+\sqrt{2+\sqrt{2}}}}}\frac{\sqrt{2-\sqrt{2+\sqrt{2+\sqrt{2}}}}}{\sqrt{2+\sqrt{2+\sqrt{2+\sqrt{2}}}}}} \\ 
 & =\text{0}\text{.82067879082866033097}\ldots \,\,,  
\end{aligned}
\]
respectively.

From these examples one can see that the set \eqref{eq_11} of the Vi\`ete-like recurrence relations gradually builds the continued fractions in the numerator and denominator of the variable ${{c}_{k}}$ at each successive step in increment of the index $k$. It is also interesting to note that each value of the variable ${{c}_{k}}$ is based on nested radicals consisting of square roots of twos only.

Figure 1 shows the dependence of the variables ${{a}_{k}}$, ${{b}_{k}}$ and ${{c}_{k}}$ as a function of the index $k$ by blue, green and red colors, respectively. We can observe how the variable ${{c}_{k}}$ tends to $1$ while the variables ${{a}_{k}}$ and ${{b}_{k}}$ tend to $2$ and $0$, respectively.
\begin{figure}[ht]
\begin{center}
\includegraphics[width=22pc]{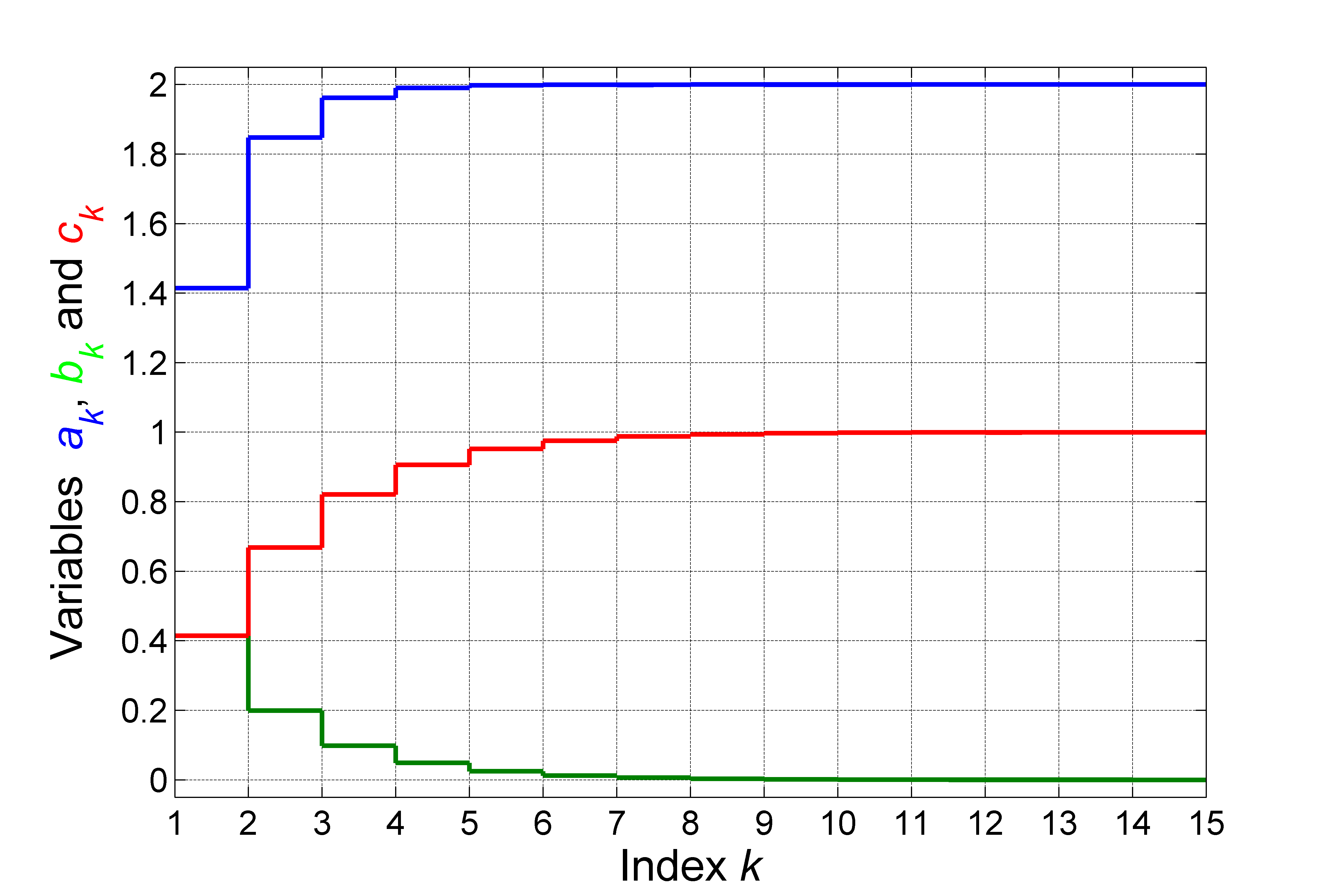}\hspace{1pc}%
\begin{minipage}[b]{28pc}
\vspace{0.3cm}
{\sffamily {\bf{Fig. 1.}} Evolution of the variables $a_k$ (blue), $b_k$ (green) and $c_k$ (red).}
\end{minipage}
\end{center}
\end{figure}

Table 1 shows the values of variable ${{c}_{k}}$ and error term ${{\varepsilon }_{k}}=1-{{c}_{k}}$ with corresponding index $k$ ranging from $4$ to $15.$ As we can see from this table, the variable ${{c}_{k}}$ quite rapidly tends to unity with increasing index $k$. In particular, the error term ${{\varepsilon }_{k}}$ decreases by factor of about $2$ at each increment of the index $k$ by one.
\begin{table}[ht]
\scriptsize
\centering
\captionsetup{width=1\textwidth}
\caption*{\small{\sffamily{\bfseries{Table 1.}} The variable $c_k$ and error term $\epsilon_k$ at index $k$ ranging from $4$ to $15$.}}
\begin{tabular}{p{0.75cm} p{4cm} p{3.5cm}}
\hline
$k$ & \qquad\qquad\,\, \ $c_k$    & \qquad\qquad\quad $\epsilon_k$ \\ [0.5ex]
\hline\hline
4   & 0.90634716901914715794...  & 0.09365283098085284205... \\
5   &	0.95207914670092534858...  & 0.04792085329907465141... \\
6   &	0.97575264993237653232...  & 0.02424735006762346767... \\
7   &	0.98780284145152917070...  & 0.01219715854847082929... \\
8   & 0.99388282491415211156...  & 0.00611717508584788843... \\
9   &	0.99693673501114949604...  & 0.00306326498885050395... \\
10  &	0.99846719455859369106...  & 0.00153280544140630893... \\
11  &	0.99923330359286120490...  & 0.00076669640713879509... \\
12  &	0.99961657831851611515...  & 0.00038342168148388484... \\
13  &	0.99980827078273533526...  & 0.00019172921726466473... \\
14  &	0.99990413079635610519...  & 0.00009586920364389480... \\
15  &	0.99995206424931502866...  & 0.00004793575068497133... \\
\hline
\end{tabular}
\normalsize
\end{table}

\section{New formula for pi}

As the error term ${{\varepsilon }_{k}}$ decreases successively by factor of about $2$ (see third column in the Table 1), we may expect that $2^k\varepsilon_k$ is convergent and tends to some constant when the index $k$ tends to infinity. The computational test shows that the value $2^k\varepsilon_k$ approaches to $\pi/2$ as the index $k$ increases. Therefore, we assume that
\[
\lim_{k \to \infty} 2^{k} \varepsilon_k = \frac{\pi}{2}
\]
or
\[
 \pi = \lim_{k \to \infty} 2^{k+1} \left(1-c_k\right).
\]
Furthermore, relying on numerical results we also suggest a generalization to the power $m$ as given by
\begin{equation}\label{eq_12}
m\,\pi = \lim_{k \to \infty} 2^{k+1} \left(1-c_k^m\right).
\end{equation}

Since the variable $c_k$ is determined within the set \eqref{eq_11} of the Vi\`ete-like recurrence relations, the new equation \eqref{eq_12} can also be regarded as the Vi\`ete-like formula for pi.

\section{Conclusion}

We show a set \eqref{eq_11} of the Vi\`ete-like recurrence relations for the constant $1$ derived by using the Vi\`ete-like formula \eqref{eq_2} for pi. Sample computations reveal that the variable ${{c}_{k}}$ quite rapidly tends to unity as the index $k$ increases.

\section*{Acknowledgments}

This work is supported by National Research Council Canada, Thoth Technology Inc. and York University.

\bigskip


\end{document}